\documentclass[11pt]{article}

\usepackage{amsthm}
\usepackage{moreverb,url}
\usepackage{graphicx}
\usepackage{subfigure}
\usepackage{caption}
\usepackage{siunitx}
\usepackage{multirow}
\usepackage{adjustbox}
\usepackage{rotating}
\usepackage{makecell}
\usepackage{amsfonts}
\usepackage{mathtools}
\usepackage{authblk}

\graphicspath{ {./images/} }
\begin{document}

\author[1]{Antonella Galizia$^{\star}$}

\author[1]{Simone Cammarasana$^{\star}$}
\author[1]{Andrea Clematis}
\author[1]{Giuseppe Patan\'e}

\affil{CNR, IMATI, Italy, Genova
 \newline *Joint first authors, with equal contribution
 }

\title{Evaluating Accuracy and Efficiency of HPC Solvers for Sparse Linear Systems with Applications to PDEs}

\maketitle

\begin{abstract}{Partial Differential Equations (PDEs) describe several problems relevant to many fields of applied sciences, and their discrete counterparts typically involve the solution of sparse linear systems. In this context, we focus on the analysis of the computational aspects related to the solution of large and sparse linear systems with HPC solvers, by considering the performances of direct and iterative solvers in terms of computational efficiency, scalability, and numerical accuracy.  Our aim is to identify the main criteria to support application-domain specialists in the selection of the most suitable solvers, according to the application requirements and available resources. To this end, we discuss how the numerical solver is affected by the regular/irregular discretisation of the input domain, the discretisation of the input PDE with piecewise linear or polynomial basis functions, which generally result in a higher/lower sparsity of the coefficient matrix, and the choice of different initial conditions, which are associated with linear systems with multiple right-hand side terms. Finally, our analysis is independent of the characteristics of the underlying computational architectures, and provides a methodological approach that can be applied to different classes of PDEs or with approximation problems.}
\end{abstract}

\section{Introduction\label{sec:INTRODUCTION}}
\begin{figure}[t] 
\centering
\begin{tabular}{ccc}
(a)\includegraphics[height=90pt]{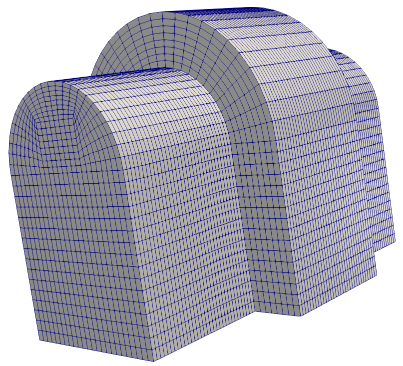}
&(b)\includegraphics[height=90pt]{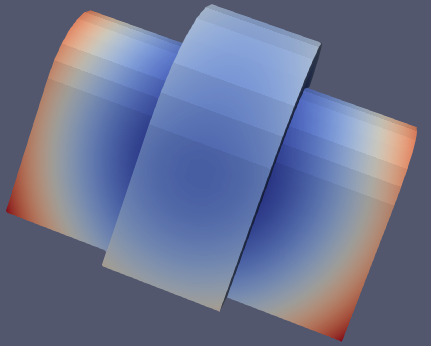}
&(c)\includegraphics[height=90pt]{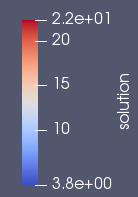}
\end{tabular}
\begin{tabular}{cc}
(d)\includegraphics[height=100pt]{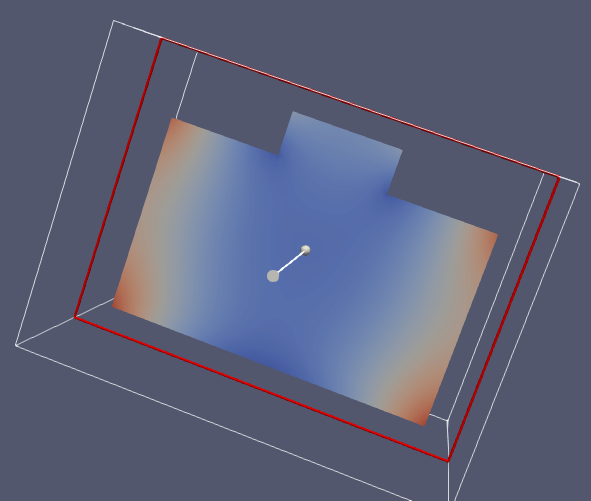}
&(e)\includegraphics[height=100pt]{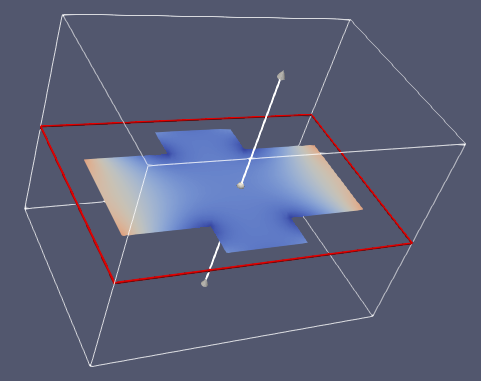}
\end{tabular}
\caption{(a) Input CAD model, (b) color-map and (c) values of the solution to the Laplace equation. (d,e) Solution on two orthogonal planes in the inner volume.\label{fig:citrusSolve}}
\end{figure}

Partial differential equations (PDEs) are crucial for the understanding of the behaviour of several phenomena~\cite{evans10} in engineering~\cite{donatelli2015robust}, geophysical exploration~\cite{keller2017performance}, and fluido-dynamics~\cite{turek2010feast}. Since it is not always possible to compute the analytic solution of PDEs, numerical approaches discretise PDEs and generally require the solution to linear systems~\cite{quarteroni2017numerical}. In this context, previous work typically applies a workflow that includes the following steps: modelling of the input problem, discretisation of the input domain and of the PDE, and solution of the discrete equation, which typically results in a large, sparse linear system (Fig.~\ref{fig:citrusSolve}). 

With respect to this framework, we focus on the analysis of computational aspects related to the solution of large and sparse linear systems with HPC solvers. We analyse the performance of direct and iterative solvers by taking into account the impact of regular/irregular discretisation grids of the input domain on the computational cost and on the numerical accuracy of the solution. As additional elements of our analysis, we discuss the influence of the discretisation of the input PDE with piecewise linear or polynomial basis functions and with different initial conditions. In fact, the discretisation of PDEs with Finite Element Methods (FEMs) and polynomial basis functions generally results in a reduced sparsity of the coefficient matrix with respect to the linear basis, and different initial conditions are associated with linear systems with multiple right-hand side terms. As test cases, we consider the Laplace equation on a 2D/3D domain, as it is a standard problem for testing numerical solvers, and it is related to several physical applications (e.g., electrostatic potential, heat diffusion in steady state). However, a similar analysis can be extended to different classes of PDEs, such as Laplacian eigen problems, heat equation, etc.

To solve sparse linear systems associated with PDEs, we employ parallel software freely accessible (SuperLU, PETSc~\cite{lidemmel03,petsc-user-ref}) and consider several configurations, performance metrics, and properties, such as efficiency, scalability, and accuracy. Instead of evaluating the impact of the underlying architecture on the performance of the solvers, our objective is to identify the main criteria to support application specialists in the selection of the most suitable solvers, according to the underlying application and available computational resources.  Following the approach proposed in~\cite{asanovic2009view}, our aim is to derive a ``rule of thumbs'' for \emph{production-oriented programmers}, i.e., domain specialists exploiting, at their best, the tools provided by \emph{optimisation-oriented programmers}, who are committed towards the design of efficient and optimised libraries. With this goal in mind, we characterise sparse linear solvers by taking into account application constraints, such as the discretisation of the input domain (e.g., regular/irregular), the properties of the coefficient matrix (e.g., size, sparsity, conditioning) and of the solver (e.g., accuracy, scalability). The proposed analysis will support the exploitation of parallel tools, as available off-the-self solvers, and possibly reduce the barrier to entry in HPC ecosystems.  Furthermore, thanks to its generality, our discussion represents a methodological approach to be applied with different classes of PDEs (e.g., elliptic, parabolic PDEs) or with approximation schemes.

\paragraph*{Related work, novelty, and main contributions\label{sec:RELATED-WORK}}
The scientific community has a great interest in the numerical solution of linear systems, as confirmed by a long-standing tradition in software development~\cite{dongarra1995software},~\cite{buttari2006impact},~\cite{song2015scalable}. An overview on the state-of-the-art of  widespread tools is presented in Sect.~\ref{sec:Solv}.

Different works propose efficient parallel solvers of sparse linear systems: novel strategies and related performance figures are compared with well-known software packages, together with the discussion of theoretical complexity and accuracy, and/or the use of specific parallel paradigms, run-time systems and accelerators. In~\cite{chen2018distributed}, a parallel hierarchical solver is proposed and compared with the SuperLU performances. In~\cite{wang2016parallel}, the performances of a novel multi-frontal solver are discussed and compared with MUMPS~\cite{MUMPS:1} and SuperLU. In~\cite{agullo2013multifrontal,agullo2016implementing}, a DAG-based parallelisation of the QR factorization of sparse matrices is implemented over the StarPU run-time system. The GHOST toolkit~\cite{kreutzer2017ghost} implements a collection of building blocks for sparse linear algebra on heterogeneous systems, such as multicore processors, graphics processing units (GPUs), and other accelerators (e.g., Intel Xeon Phi). Focusing on GPUs, in~\cite{naumov2011incomplete}, NVIDIA GPU-accelerated libraries are used to implement the incomplete-LU and Cholesky preconditioned iterative solvers; in~\cite{abdelfattah2016performance} and~\cite{liu2017fast}, optimised GPU implementations of key operations in sparse linear systems solvers are proposed. In~\cite{dongarra2016high}, a new benchmark for the computation of the \emph{High-Performance Conjugate Gradient} for the solution of the Poisson equation on a regular 3D grid is presented. This work also discusses a closer relation of the benchmark with application requirements to drive future systems' design.

This previous work emphasises the need of effective sparse linear solvers from the developer's viewpoint, i.e., optimisation-oriented programmers investigating the development of optimised solvers with notable performances, the exploitation of accelerators, and the design of future systems. From a different perspective, our aim is to address the production-oriented programmers needs through a comparison of well-known tools and the identification of the main criteria for the selection of the direct or iterative solvers of sparse linear systems that better suit application requirements.

More close to our approach, in~\cite{puzyrev2016evaluation, carracciuolo2011computational} several parallel solvers of sparse linear systems discretising PDEs are applied for applications in geophysical exploration and fluid-dynamics, also discussing solvers' effectiveness with respect to the specific applications. In~\cite{carracciuolo2011computational}, a parallel software environment for 3D fluid dynamic simulations is presented; three case studies are proposed to test direct and iterative solvers on different levels of simulation complexity. While these works discuss solvers for specific problems, we aim at addressing the needs of a larger class of PDEs. To the best of our knowledge, an up-to-date comparison of the computational efficiency and numerical accuracy of direct and iterative solvers of large sparse linear systems (e.g., associated with discrete PDEs) is actually missing. 
\paragraph*{Paper organisation}
Firstly (Sect.~\ref{sec:EXP-PIPELINE1}), we present the selected case study, the parallel tools employed to solve the corresponding sparse linear systems, data sets with different matrix sparsity, metrics and test-beds. Then, we present the performance figures of direct, and iterative solvers (Sect.~\ref{sec:EXPERIMENT-RESULTS}), which provide a fine granularity of analysis about scalability, and possible bottlenecks of each solver. Solvers are then evaluated in terms of memory requirements, numerical accuracy of the solution, multiple right-hand side vectors and linear/polynomial basis functions (Sect.~\ref{sec:COMPARISON}), also presenting their main limitations.  Finally, (Sect.~\ref{sec:DISCUSSION}) we compare numerical solvers and we present the rule of thumbs to support application specialists; last considerations and future work conclude the proposed analysis (Sect.~\ref{sec:CONCLUSION}).

\section{Experimental pipeline}\label{sec:EXP-PIPELINE1}
We briefly review the main aspects related to direct and iterative solvers (Sect.~\ref{sec:Solv}), our case studies (Sect.~\ref{sec:CASE-STUDIES}), metrics and testbeds (Sect.~\ref{sec:metrics}).

\subsection{Direct and iterative solvers\label{sec:Solv}}
\paragraph*{Numerical libraries}
Numerical linear algebra has always represented a central topic for the scientific community; many high-performance software packages and tools for parallel heterogeneous architectures have been developed in the last decades~\cite{dongarra1995software},~\cite{buttari2006impact},~\cite{dongarra2016parallel}. Several studies propose general purpose solvers and preconditioners of sparse linear systems on parallel resources. For instance, PaStiX~\cite{henon2002pastix}, SuperLU~\cite{li2005overview}, MUMPS~\cite{MUMPS:1}, and PARDISO~\cite{schenk2004solving} provide parallel sparse direct solvers. Hypre~\cite{falgout2006design} is a parallel suite of preconditioners and solvers for sparse linear systems. Trilinos~\cite{heroux2005overview} and PETSc  - Portable, Extensible Toolkit for Scientific Computation~\cite{petsc-user-ref} represent frameworks for the solution of complex scientific problems; among others, both libraries provide several modules and external packages for sparse iterative and direct solvers. 

Since all the aforementioned software tools have been extensively used, updated, and tested on a large set of problems, they represent the ideal base for our discussion, as we are interested in a systematic analysis of the performance and accuracy of available parallel solvers for sparse linear systems. To this end, we experiment different algorithms and tools~\cite{cammarasanahigh}, compare MUMPS~\cite{MUMPS:1,MUMPS:2} and SuperLU~\cite{lidemmel03} for direct solvers, and focus on PETSc~\cite{petsc-user-ref} for iterative solvers and preconditioners. 

\paragraph*{Direct and iterative linear solvers}
To discuss performance and accuracy figures, we briefly review the main aspects and features of direct and iterative solvers of a linear system \mbox{$\textbf{A}\textbf{u}=\textbf{b}$}. \emph{Direct solvers}~\cite{golub2012matrix} are based on the LU \emph{decomposition}, which factorises the coefficient matrix as~$\textbf{A}=\textbf{L}\textbf{U}$, where~$\textbf{L}$ and~$\textbf{U}$ are a lower and an upper triangular matrix, respectively. The triangular linear systems~$\textbf{L}\textbf{y}=\textbf{b}$,~$\textbf{U}\textbf{u}=\textbf{y}$ are solved with a forward and a backward substitution in linear time with respect to the number of unknowns, thus reducing the quadratic or cubic computational cost for the solution to the input linear system.

If~$\textbf{A}$ is sparse, then the factorisation typically introduces many more nonzero entries, i.e., the \mbox{$(i,j)$} entry of~$\mathbf{A}$ can be non-zero in the~$(\mathbf{L},\mathbf{U})$ factors even if it is originally zero. This phenomenon is known as \emph{fill-in} and has an impact on (i) the growth of the memory requirement and the number of operations, (ii) an a-priori unknown structure of the factorised matrices. Re-ordering methods (e.g., approximate minimum degree method~\cite{amestoy1996approximate}, multilevel graph partitioning~\cite{karypis1998fast}) can be applied to the columns of the coefficient matrix to minimise the fill-in introduced by the LU factorisation. Moreover, a symbolic factorisation can be used to define a sparsity pattern of the LU decomposition with a reduced computational complexity with respect to the numerical factorisation. To efficiently generate the ($\textbf{L}$,~$\textbf{U}$) matrices, we can apply the left-looking (Spooles~\cite{ashcraft2002solving}), right-looking (SuperLU\_DIST), and multi-frontal (MUMPS) method. 

\emph{Iterative solvers}~\cite{golub2012matrix,saad2003iterative} (e.g., Krylov subspace projection, Arnoldi or Lanczos iterations, transpose-free variants) compute an approximated solution of a linear system, starting from an initial guess and without factorising the input coefficient matrix. The efficiency of a solver depends mainly on the number of operations for each iteration and on the number of iterations needed to converge to the solution under a certain stopping criterion. The number of iterations depends also on the condition number of the coefficient matrix, which can be reduced with a preconditioning operation. We briefly recall that preconditioners are matrices used to rewrite the input linear system in an equivalent formulation, whose new coefficient matrix has a lower conditioning. Given a linear system \mbox{$\textbf{A}\textbf{u}=\textbf{b}$} and a preconditioner~$\mathbf{M}^{-1}$, we solve the linear systems
\begin{equation*}
\textbf{A}\mathbf{M}^{-1}\textbf{y}=\textbf{b},\qquad
\textbf{u}=\mathbf{M}^{-1}\textbf{y},\qquad
\kappa(\textbf{A}\mathbf{M}^{-1})<<\kappa(\textbf{A}).
\end{equation*}
Besides the computation of the matrix~$\mathbf{M}^{-1}$, preconditioners generally improve the numerical stability of direct and iterative solvers; for iterative solvers, preconditioners generally increase the number of operations performed for each iteration but reduce the number of iterations necessary to converge.

\subsection{Case studies}\label{sec:CASE-STUDIES}
\begin{figure}[t]
\centering
\begin{tabular}{cc}
(a)\includegraphics[height=130pt]{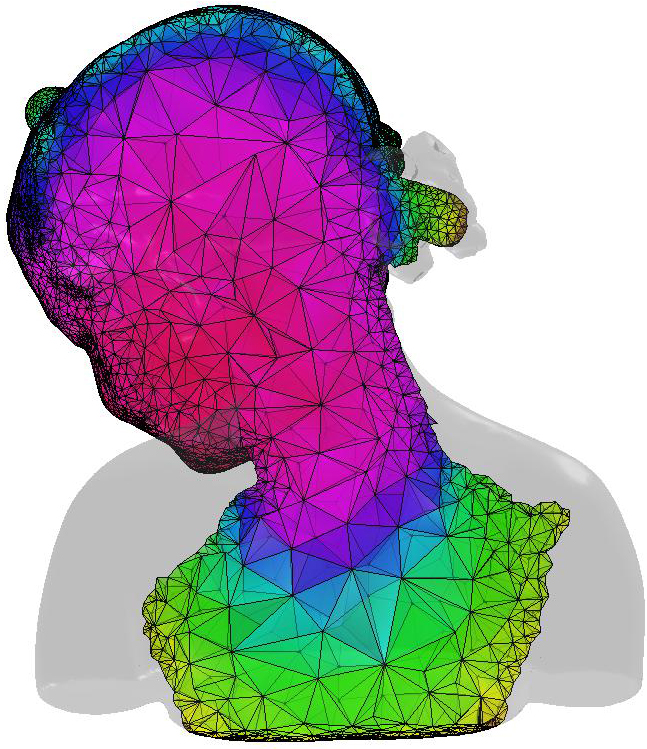}
&(b)\includegraphics[height=130pt]{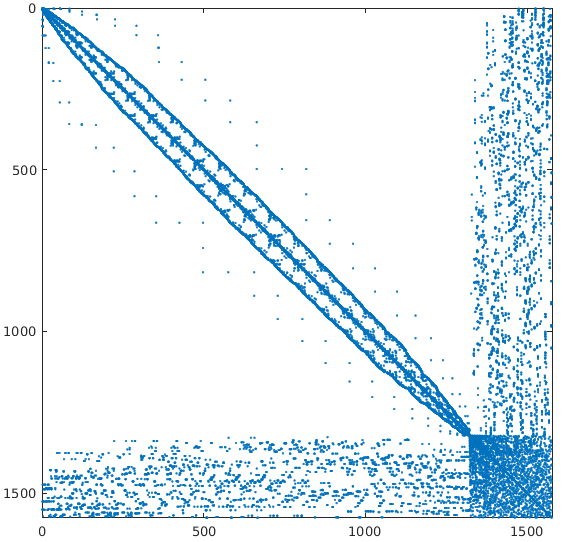}
\end{tabular}
\caption{(a) Tetrahedralisation of a 3D domain and (b) sparsity pattern of the corresponding Laplacian matrix.\label{fig:MESH-SPARSITY}}
\end{figure}

As case study, we consider the Laplace equation~$\Delta u=0$ on a 2D/3D domain~$\Omega$ with Dirichlet condition~$u=f$ on the boundary of~$\Omega$. If the input discrete domain is represented as a regular grid, then the Laplace equation is discretised with a finite difference scheme and the corresponding coefficient matrix is sparse and has a regular structure. If the input domain is discretised as an irregular grid (e.g., a tetrahedral mesh as in Fig.~\ref{fig:MESH-SPARSITY}(a)), then we apply a finite element discretisation~\cite{quarteroni2017numerical}; the corresponding coefficient matrix is still sparse but it has an irregular structure. For tetrahedral meshes, the Laplacian matrix \mbox{$\mathbf{A}:=\mathbf{B}^{-1}\mathbf{L}$} is defined as the product between the diagonal matrix~$\mathbf{B}$, whose entries are the sum of the volume of the tetrahedra incident at a given vertex, and~$\mathbf{L}$ is the Laplacian matrix with entries \mbox{$L(i,j):=w(i,j):=\frac{1}{6}\sum_{k=1}^{n}l_{k}\cot\alpha_{k}$} for each edge \mbox{$(i,j)$}, \mbox{$L(i,i):=-\sum_{j\in N(i)}w(i,j)$}, and zero otherwise~\cite{patane2016star}. For regular grids,~$w(i,j):=1$. In all the aforementioned cases, the Laplace-Beltrami operator is discretised as the~$n\times n$ \emph{Laplacian matrix}
\begin{equation*}
A(i,j):=
\left\{
\begin{array}{ll}
w(i,j) &(i,j) \textrm{ edge};\\  
-\sum_{k}w(i,k) &i=j,
\end{array}
\right.
\end{equation*}
which is sparse and positive semi-definite, but not necessarily symmetric. Then, the Laplace equation is discretised as the sparse linear system~$\textbf{A}\textbf{u}=\textbf{b}$, where~$\textbf{u}$ is the array of the values of the solution at the nodes of the input grid and~$\mathbf{b}$ is the right-hand side term that defines the boundary condition. Finally, we do not make assumptions on the coefficient matrix and/or constraints on the connectivity of the input grid, such as a maximum degree of the nodes. 

Recalling that the \emph{non-zero density} of a matrix is defined as the percentage of its non-zero entries with respect to the total number of elements of the matrix and the \emph{sparsity pattern} as the structure of the non-zero entries in the coefficient matrix, the Laplacian matrix on irregular grids (Table~\ref{tab:DomainandGrid}) has an arbitrary sparsity pattern and a non-zero density percentage that is generally higher than on regular grids, as a matter of the generally higher number of edges incident to a given node. Indeed, the sparsity pattern of the Laplacian matrix associated with irregular grids has an impact on the data distribution among processes, MPI communications, and processes work load. However, the reduced non-zero density percentage can be exploited to improve the arithmetic density of the operations through an appropriate column reordering. These considerations further motivate a comparison of solvers of sparse linear systems associated with PDEs discretised on a regular/irregular grid.
\begin{table}[t]
\centering
\caption{Grid size and sparsity of the coefficient matrix of the input sparse linear systems.\label{tab:DomainandGrid}}
{
\resizebox{\textwidth}{!}{
\begin{tabular}{|l|l|l|l|l|}
\hline
\multicolumn{5}{|c|}{\textbf{ Regular grids }}\\ \hline
\textbf{Domain} &$\sharp$\textbf{Nodes} &\textbf{$\sharp$ Mat. Rows} &\textbf{$\sharp$ Non-zeros} &\textbf{Non-zeros} \% \\ \hline
Cube~1 &~$128 \times 128 \times 128$  &2\,097\,152 &14\,099\,408 &~$\num{3e-4}$ \\ \hline
Cube 2 &~$256 \times 256 \times 256$ & 16\,777\,216 &115\,099\,600 &~$\num{4.1e-5}$ \\ \hline
Cube 3 &~$512 \times 512 \times 512$ &134\,217\,728 &930\,123\,728 &~$\num{5e-6}$\\ \hline \hline
\multicolumn{5}{|c|}{\textbf{ Irregular grids }}\\ \hline 
\textbf{Domain} &$\sharp$\textbf{Nodes} &\textbf{$\sharp$ Mat. Rows} &\textbf{$\sharp$ Non-zeros} &\textbf{Non-zeros} \% \\ \hline
Sphere~1 &2\,094\,977 &2\,094\,977 &33\,225\,967 &$\num{7e-4}$ \\ \hline
Sphere 2 &2\,094\,834 &2\,094\,834 &61\,055\,286 &$\num{1.4e-3}$ \\ \hline
\end{tabular}}
}
\end{table}
\subsection{Metrics and testbeds\label{sec:metrics}}
To evaluate the solvers' performances, we measure the following metrics
\begin{itemize}
\item \emph{execution time} expressed in seconds;
\item \emph{FLOPS} - FLoating point Operations Per Second - expressed in Giga;
\item \emph{efficiency$(n)$}=$\frac{speed up(n)}{n}$ with~$n$ processes;
\item \emph{FLOPS efficiency$(n)$}$=\frac{F(n)}{F(1)}/n$ where~$F(1)$ and~$F(n)$ are FLOPS with~$1$ and~$n$ processes, respectively; 
\end{itemize}
In order to discuss the accuracy of direct and iterative solvers, we measure the approximation error between the exact and the computed solution as
\begin{equation}\label{eq:ERROR-METRIC}
x_{err}
=\frac{\left\|\mathbf{x}_{ground\_truth}-\mathbf{x}_{computed}\right \|_2}{\left\| \mathbf{x}_{ground\_truth}\right \|_2},
\end{equation}
where the discrete ground-truth solution of the input PDE is given by sampling the analytic solution at the vertices of the input domain. As exit condition of iterative solvers, we consider the relative error
\begin{equation}\label{eq:relativeError}
\frac{\left \|\mathbf{A}\mathbf{x}-\mathbf{b}\right \|_2 }{\left \|\mathbf{b}\right \|_2}\leq\epsilon,
\end{equation}
where~$\epsilon$ is a given threshold.

The memory requirement of each method is measured by the total memory (i.e., including support matrices and vectors) allocated in the RAM. The allocated memory is expressed in MB. These metrics will be used for an analysis of the performances (Sect.~\ref{sec:EXPERIMENT-RESULTS}) and for a comparison (Sect.~\ref{sec:COMPARISON}) of direct and iterative solvers, in order to identify which is more suitable on the basis of the application and the resources available. All performance tests have been run on the CINECA cluster Marconi, based on Intel Xeon product family and classified at the 12-th position on Top500 ranking~\cite{strohmaier2006top500} in November 2016. In particular, we have considered the Broadwell partition, composed by 2x18-cores Intel Xeon E5-2697 v4 (Broadwell) at 2.30 GHz, 1512 nodes, 36 cores/node. Total cores: 54432; 128 GB/node of RAM. 
\begin{figure}[t]
\centering
\includegraphics[height=190pt]{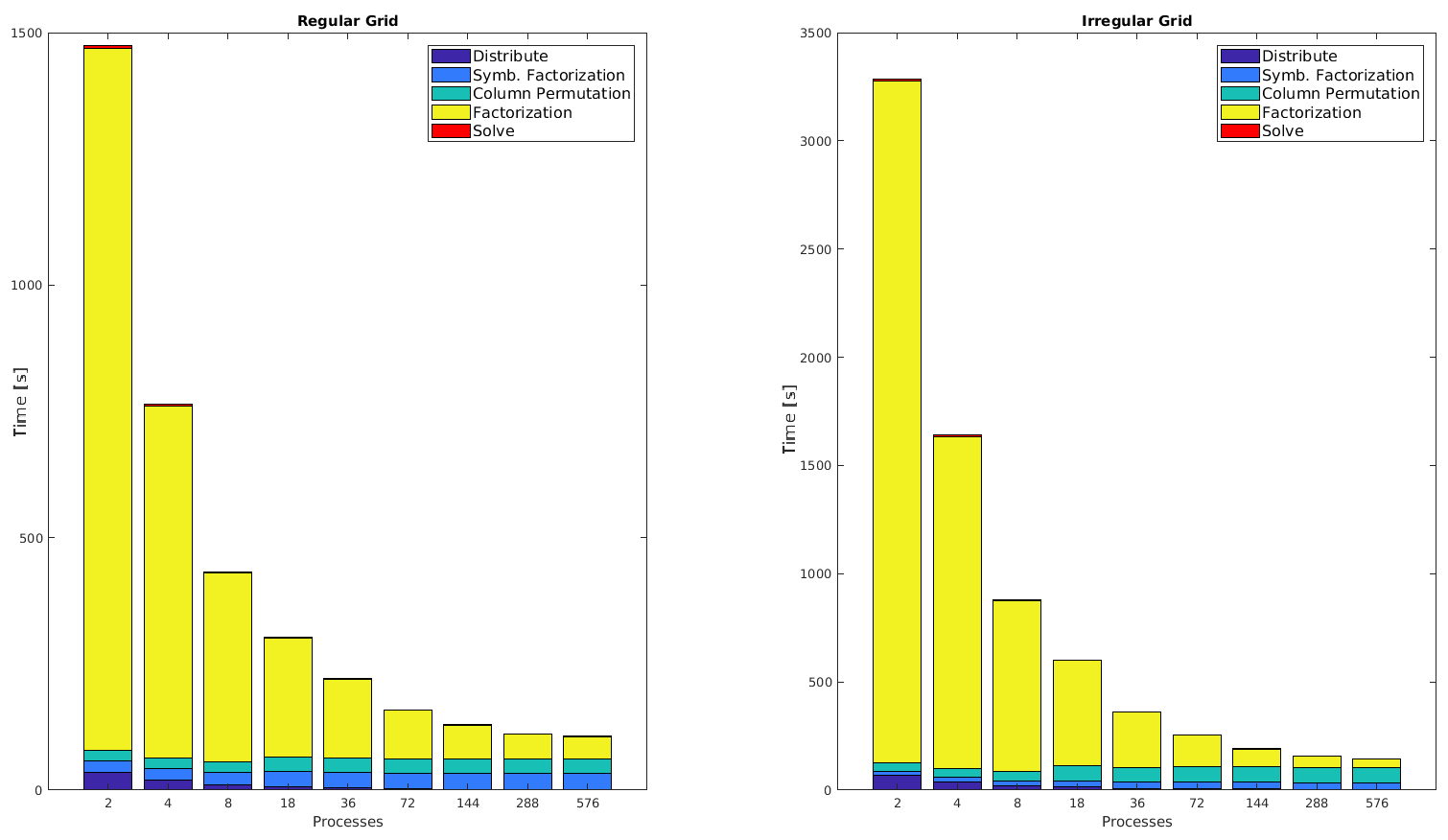}
\caption{SuperLU execution time profiled according to its five main operations.\label{fig:DM_Prof}}
\end{figure}
\section{Experimental results}\label{sec:EXPERIMENT-RESULTS}
We discuss our experiments on direct (Sect.~\ref{sec:DIRECT-METHODS}) and iterative (Sect.~\ref{sec:ITERATIVE-METHODS}) solvers.

\subsection{Experimenting direct solvers\label{sec:DIRECT-METHODS}}
According to our preliminary comparison~\cite{cammarasanahigh} between MUMPS~\cite{MUMPS:1, MUMPS:2} and SuperLU\-\_DIST\footnote{Version SuperLU\_DIST 5.2.0} is the parallel extension to the serial SuperLU for parallel distributed memory architectures~\cite{lidemmel03}. The latter obtained better performance and has been selected for a further investigation in this work, and for the sake of simplicity and readability since now on we will refer it as SuperLU. It subdivides the computation in five main steps
\begin{enumerate}
\item \emph{data distribution} among processes, which implements a cyclic block pattern of different size, according to the non-zero patterns of the input coefficient matrix;
\item \emph{column permutation}, which reorders the columns of the coefficient matrix to optimise the number of non-zero entries after factorisation;
\item \emph{symbolic factorisation}, which provides the number and position of non-zero entries of the ($\mathbf{L},\mathbf{U}$) matrices;
\item \emph{factorisation}, which factorises the coefficient matrix in the LU form, where~$\mathbf{L}$ and~$\mathbf{U}$ are a lower and upper triangular matrix, respectively;
\item \emph{solver}, which computes the solution of the triangular systems.
\end{enumerate}
As input data sets (Table~\ref{tab:DomainandGrid}), we consider a regularly tessellated cube and an irregularly tessellated sphere. According to Table~\ref{tab:dom&fact}, the LU factorisation increases the fill-in ratio of the input coefficient matrix by 2 and 3 orders of magnitude for regular and irregular domains, respectively. This phenomenon leads to higher memory and computational costs for the solution of the corresponding linear systems.
\begin{table}[t]
\centering
\caption{Size, non-zeros density and fill-in percentage of the coefficient matrix and its LU factorisation.\label{tab:dom&fact}}
{\begin{tabular}{|c|c|c|c|}
\hline
\textbf{Domain} &\textbf{A rows} & \textbf{A Non-zeros} \% &\textbf{LU fill-in \%}\\ \hline
Cube~1 & 2097152 & ~$\num{3.2e-6}$ & ~$\num{9.6e-4}$ \\ \hline 
Sphere~1 & 2094977 &~$\num{7.5e-6}$ &~$\num{1.7e-3}$ \\ \hline
\end{tabular}}
\end{table}

In Fig.~\ref{fig:DM_Prof}, we report the execution time achieved by varying the number of processes up to 576 processes and profiled on the base of the aforementioned operations. Due to memory requirements for the allocation of the factorised matrices, it was not possible to solve the system with a single processor; indeed, for the performance evaluation we compare the execution time on~$p$ processors with respect to the execution time on~$p_1$ processors,~$p_1>1$. Memory allocation for both direct and iterative solvers are discussed in Sect.~\ref{sec:grids-comparison}. 

The system resolution on a regular grid requires more than 20 minutes with 2 processes, reduces to less than 1 and a half minute with 72 processes, and 45 seconds with 576 processes. For irregular grids, the numerical solver requires more than 50 minutes with 2 processes, reduces to 2 minutes and a half with 72 processes, and 37 seconds with 576 processes. Growing from 288 to 576 processes has a little impact on time reduction, especially for regular grids. This aspect is partially due to the sequential executions of the \emph{column permutation} and of the \emph{symbolic factorisation}. Since these operations do not scale on the processes, their impact on the global computation becomes heavier while increasing the number of processes. On regular grids, the aggregated impact of both operations is 3\% with 2 processes, and it increases to 55\% with 576 processes; on irregular grids, it grows from 1.7\% with 2 processes to 62\% with 576 processes. On the contrary, \emph{data distribution} and \emph{solver} smoothly scale and have a negligible effect. In fact, the aggregated impact counts less than 3\% on the global time. Indeed, the \emph{factorisation} dominates algorithm' complexity,~\cite{lidemmel03,li2005overview,MUMPS:3}. 

To further investigate the \emph{factorisation} operation (yellow parts of Fig.~\ref{fig:DM_Prof}), in Table~\ref{tab:DMFact} we report FLOPS and efficiency on regular and irregular grids obtained by varying the number of processes. We note that the metric is normalised with 2 processes, as memory allocation problems occurred with 1 process. The values increase from about 32 GFLOPS on 2 processes to about 1 TFLOPS (11.1\% efficiency) and to 2.7 TFLOPS on 576 processes (30\% efficiency) on regular and irregular domains, respectively. It follows that the FLOPS metrics is interesting even if the efficiency is limited, since reaching the Tera scale is a plus. In general, better results are achieved on irregularly tessellated domains, as they have an improved arithmetic density for direct solvers, while maintaining limited the impact of the memory pattern and communication overhead.
\begin{table}[t]
\centering
\caption{Direct solvers, FLOPS, and factorisation efficiency.\label{tab:DMFact}}
{
\begin{tabular}{|c|cc|cc|}
\hline
\textbf{Processes} & \multicolumn{2}{c|}{\textbf{Regular grid}} & \multicolumn{2}{c|}{\textbf{Irregular grid}} \\ 
\hline
{} & GFLOPS & Efficiency [\%] & GFLOPS & Efficiency [\%] \\ \hline
1 & - & - & - & - \\ \hline
2 &  32 & 100 & 32 & 100 \\ \hline
4 &  64 & 99.8 & 66 & 100  \\ \hline
8 &  119 & 93.2 & 129 & 99.6 \\ \hline
18 &  189 & 65.6 & 272 & 92.8 \\ \hline
36 &  288 & 49.8 & 504 & 85.9  \\ \hline
72 & 467 & 40.4 & 758 & 64.6 \\ \hline
144 &  670 & 29.0 & 1330 & 56.6 \\ \hline
288 &  903 & 19.5 & 2073 & 44.2 \\ \hline
576 &  1033 & 11.1 & 2765 & 29.4 \\ \hline
\end{tabular}}
\end{table}
\subsection{Experimenting iterative solvers\label{sec:ITERATIVE-METHODS}}
\begin{figure}[t]
\centering
\includegraphics[width=0.9\textwidth]{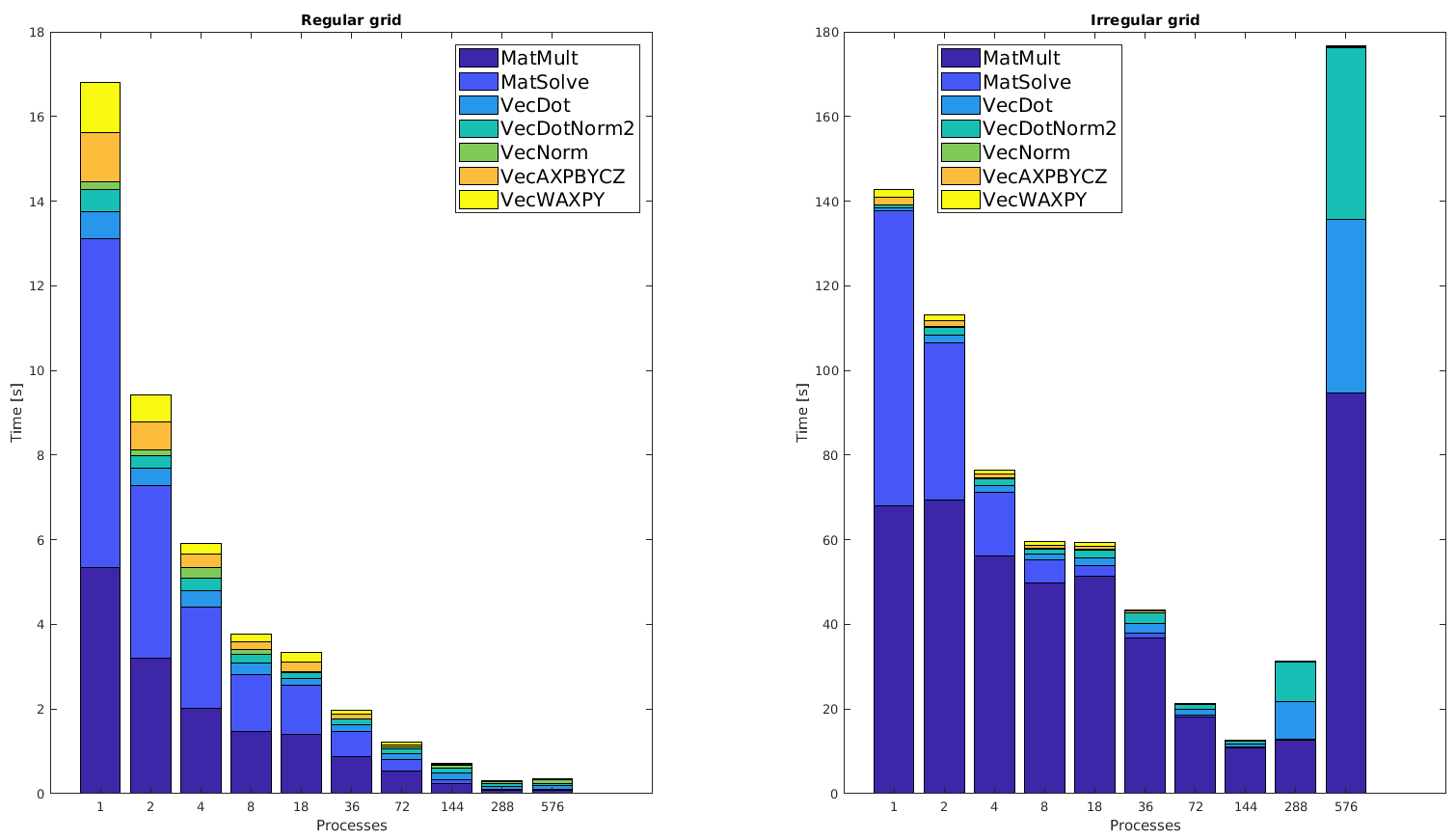}
\caption{BICGSTAB execution time profiled with matrix and vector operations.\label{fig:IterativeProfiling}} 
\end{figure}

We selected PETSc~\cite{petsc-user-ref} to test 5 iterative solvers (BICGSTAB, GMRES, IBICGSTAB, TFQMR, CGR) of sparse linear systems without applying  preconditioners and combined with 3 preconditioners (Hypre, ASM, Block-Jacobi)~\cite{saad2003iterative}. According to the results in~\cite{cammarasanahigh}, the combination of the iterative BICGSTAB method, which is based on Lanczos bi-orthogonalisation, with the Block-Jacobi preconditioner has been selected as solver to be further investigated, as it provided better performance figures. We briefly recall that Block-Jacobi is a general purpose preconditioner, which reduces the condition number of the coefficient matrix without optimising its sparsity pattern and without affecting the non-zero density. Indeed, we test the solver on the \emph{Cube~1} and \emph{Spere 1} domains, and consider its convergence to a solution with~$\epsilon = 10^{-12}$ in Eq. (\ref{eq:relativeError}). In Fig.~\ref{fig:IterativeProfiling}, we report the execution time of the BICGSTAB algorithm profiled on the basis of the main matrix and vector operations. Execution time is very low on both regular and irregular grids, however, the algorithm has a different behaviour on regular and irregular grids. In the first case, the algorithm requires about 17 seconds to converge to a solution with 1 process, and less than half second with 576 processes; this smooth scalability can be improved on higher matrix dimension, i.e. larger data sets, as discussed in Sect.~\ref{sec:grids-comparison}. On irregular grids, BICGSTAB requires about 2 and half minutes to converge with 1 process, reduced to 13 seconds with 144 process; however, with 576 increases up to about 3 minutes. In fact, some of the main operations (e.g., \emph{Matmult} and \emph{Vecdot}, a matrix-vector multiplication and a vector scalar product, blue and light-blue parts in Fig.~\ref{fig:IterativeProfiling}) do not properly scale when increasing the number of processes.  

In Table~\ref{tab:iterativeSolve}, we report the performances of BICGSTAB on regular and irregular grids in terms of FLOPS, related efficiency, and number of iterations; on regular grids, the algorithm has an efficiency of 35\% (13.5 GFLOPS) with 36 processes and of 17\% (106 GFLOPS) with 576 processes. On irregular grids, the trend is significantly worst, with an efficiency of 17\% (2.5 GFLOPS) with 36 processes and lower than 1\% (and 0.5 GFLOPS) with 576 processes. Furthermore, the number of iterations needed to the BICGSTAB algorithm to converge grows as we increase the number of processes; this aspect has an impact on the solving time and, consequently, on the algorithm efficiency. 
\begin{table}[t]
\centering
\caption{FLOPS, efficiency, and number of iterations of BICGSTAB on regular and irregular grids.\label{tab:iterativeSolve}}
{
\resizebox{\textwidth}{!}{
\begin{tabular}{|c|ccc|ccc|}
\hline
& \multicolumn{3}{c|}{\textbf{Regular grid} - Cube~1} & \multicolumn{3}{c|}{\textbf{Irregular grid} - Sphere~1} \\ \hline
Processes & GFLOPS & Iterations & FLOPS Eff. [\%] & GFLOPS & Iterations & FLOPS Eff. \\ \hline
1  & 1.1 & 128 & 100 & 0.4 & 195 & 100  \\ \hline
2  & 2.4 & 152 & 108 & 0.8 & 275 & 99  \\ \hline
4  & 4.3 & 169 & 98 & 1.3 & 463 & 80 \\ \hline
8  & 6.3 & 168 & 72 & 1.6 & 504 & 50  \\ \hline
18 & 6.8 & 164 & 35 & 1.7 & 546 & 23  \\ \hline
36 & 13.2 & 183 & 34 & 2.5 & 574 & 17  \\ \hline
72  & 21.6 & 205 & 28 & 4.2 & 528 & 14  \\ \hline
144 & 45.2 & 242 & 29 & 8.0 & 591 & 14  \\ \hline
288  & 92.2 & 231 & 30 & 3.1 & 521 & 3  \\ \hline
576   & 104.0 & 270 & 17 & 0.5 & 539 & 1  \\ \hline
\end{tabular}}
}
\end{table}
\begin{figure}[t]
\centering
\subfigure[Regular grids]{\label{fig:matmultReg}\includegraphics[height=0.495 \linewidth]{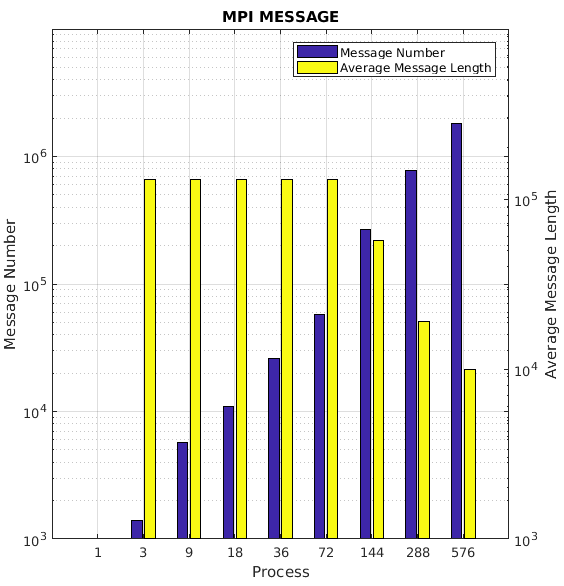}}
\subfigure[Irregular grids]{\label{fig:matmultIrreg}\includegraphics[height=0.495 \linewidth]{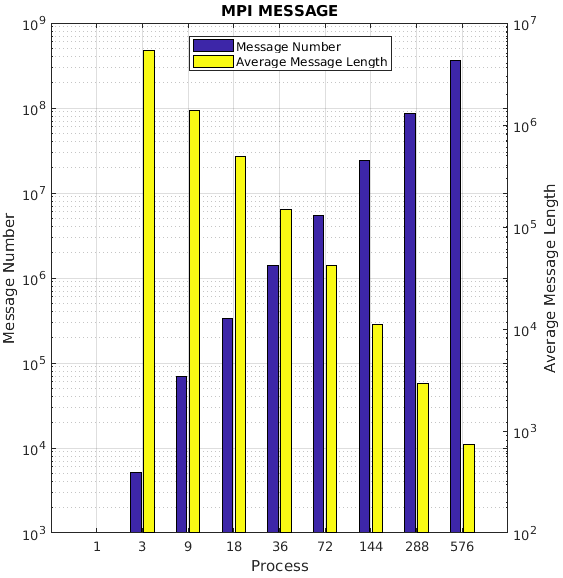}}
\caption{Matmult communications and messages' length for regular/irregular grids.}
\label{fig:matmult}
\end{figure}
Since the lower regularity of the sparsity pattern of the coefficient matrix has a deep impact on the slow convergence of the solver, we analyse the trend of the main operations. In Fig.~\ref{fig:matmult}, we propose an example of the communication figures in term of MPI message number and average message length for the \emph{Matmult} operation. Increasing the number of processes, the number of communications increases while their size reduces. However, on regular grids (a) the gap is closer than on irregular grids (b), In fact, the former message length varies from an average of~$100$ KB to~$10$ KB and the message number varies from~$10^3$ to~$10^6$, as for the latter message length varies from an average of~$1$ MB to~$1$ KB and message number varies from~$10^3$ to~$10^8$. The overhead paid to send a high number of short messages (Fig.~\ref{fig:matmultIrreg}, last column) is not balanced by a more intensive computational counterpart neither by a more accurate solution leading to a faster convergence. Actually, this trend is confirmed also on more demanding linear systems, i.e. on \emph{Sphere~2}.

\section{Computational performances and numerical accuracy}\label{sec:COMPARISON}
Besides scalability and performance figures, we propose additional experiments to characterise the properties of linear and iterative solvers peculiarities in terms of memory requirements (Sect.~\ref{sec:grids-comparison}), approximation error (Sect.~\ref{sec:ERROR-ANALYSIS}), and multiple right-hand side terms (Sect.~\ref{sec:MULTIPLE-RHS-TERMS}). We also discuss how the choice of a polynomial instead of a linear basis of the FEM method affects the performances of iterative solvers.

\subsection{Memory requirements: regular versus irregular grids\label{sec:grids-comparison}}
We analyse memory requirements of linear solvers, as they represent a bottleneck for the execution for direct solvers and for the scalability of iterative solvers.
\begin{table}[t]
\centering
\caption{Memory requirements (expressed in MB) of direct and iterative solvers on the regular grid \emph{Cube~1} and the irregular grid \emph{Sphere~1} -  average per process. \label{tab:memory regular/irregular grid}}
{
\begin{tabular}{|l|ll|ll|}
\hline
\textbf{Processes} &\multicolumn{2}{c}{\textbf{Regular grid}} & \multicolumn{2}{|c|}{\textbf{Irregular grid}}\\
\hline
{} & Direct [MB] & Iterative [MB] & Direct [MB] & Iterative[MB] \\ \hline
1 & - & 588 & - & 1020  \\ 
2 & 17976 & 340 & 37603 & 529 \\
4 & 9284 & 196 & 17512 & 298 \\
8 & 4899 & 123 & 8471 & 192 \\
18 & 2405 & 82 & 3972 & 134  \\
36 & 1393 & 66 & 1850 & 112 \\
72 & 831 & 58 & 1198 & 76 \\
144 & 580 & 53 & 694 & 61  \\
288 & 417 & 50 & 428 & 54 \\
576 & 369 & 49 & 314 & 53 \\ \hline
\end{tabular}}
\end{table}
Table~\ref{tab:memory regular/irregular grid} reports the allocated memory of direct and iterative solvers for \emph{Cube~1} and \emph{Sphere~1}, expressed as the average memory allocated per process. Both methods have a good distribution ratio of data among processes; average and maximum values are quite close, apart from few cases for direct solvers. However, direct solvers allocate one or two orders of magnitude more than iterative solvers, as a very large amount of memory is allocated for the factorised matrices; for instance, direct solvers with two processes need 17.5 GB of memory and 37 GB versus 340 MB and 530 MB required by iterative solvers on regular and irregular grids, respectively. Indeed, in the latter case both methods double the memory costs, without pushing to a higher order.
\begin{table}[t]
\centering
\caption{Scalability results of iterative solvers when increasing the dimension of the coefficient matrix, with Time, GigaFLOPS (GF.) and FLOPS efficiency (F. Eff.) metrics.\label{tab:granularity}}
{
\resizebox{\textwidth}{!}{
\begin{tabular}{|c|ccc|ccc|ccc|}
\hline
Domain &\multicolumn{3}{c|}{\textbf{Cube~1}}&\multicolumn{3}{c|}{\textbf{Cube~2}}&\multicolumn{3}{c|}{\textbf{Cube~3}}\\
\hline
Procs & Time[s] & GF. & F. Eff.[\%] & Time[s] & GF. & F. Eff.[\%] & Time[s] & GF. & F. Eff.[\%] \\ \hline
1 & 17.4 & 1 & 100 & 1047.1 & 0.3 & 100 & 4361.8 & 0.3 & 100 \\ \hline
144 & 0.73 & 46 & 29 & 12.2 & 34 & 72 & 141 & 41 & 89  \\ \hline
288 & 0.32 & 94 & 30 & 9.7 & 55 & 57 & 123 & 58 & 58 \\ \hline
576 & 0.34 & 106 & 17 & 6.7 & 81 & 42 & 73.0 & 110 & 54  \\ \hline
\end{tabular}}
}
\end{table}
From these tests, it follows that iterative solvers definitely overwhelm direct solvers in terms of memory cost, and represent the unique possible choice to manage large coefficient matrices. In fact (Sect.~\ref{sec:DIRECT-METHODS}), direct solvers have incurred in memory limits for the computation of the LU factorisation when largely increasing the number of nodes of the input grid. On the contrary, iterative solvers generally do not suffer of memory limits, once the coefficient matrix has been allocated.  

\begin{figure}[t]
\centering
\includegraphics[height=180pt]{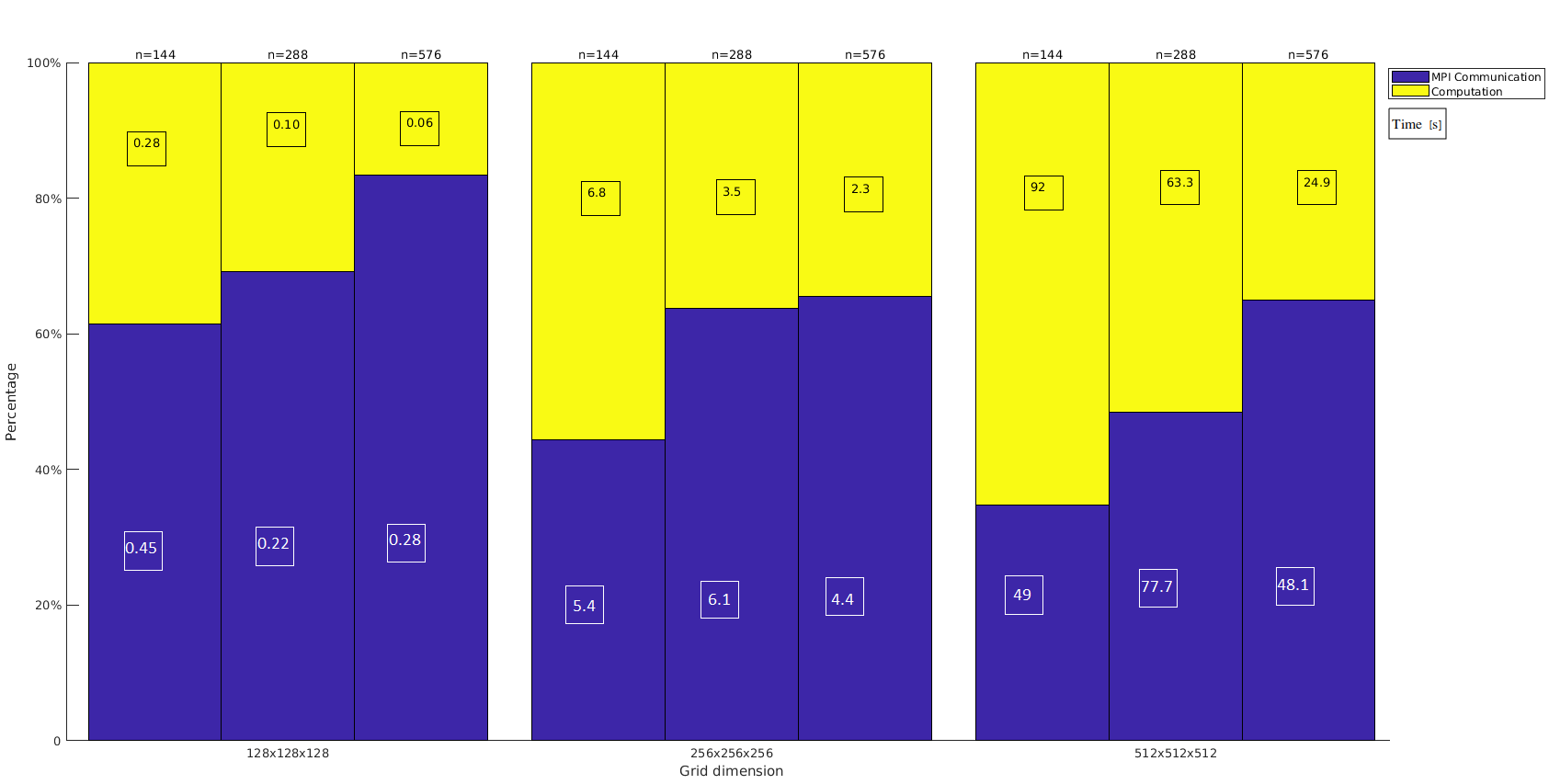}
\caption{Timings of the computation and communication phases by varying the number of processes and the size of the coefficient matrix; tests consider regular grids.\label{fig:granularity}} 
\end{figure}

We further test iterative solvers considering larger regular domains. According to the results in Table~\ref{tab:granularity}, the scalability of the algorithm improves with the dimension of the coefficient matrix: BICGSTAB requires (i) 18 minutes with 1 process and 7 seconds with 576 processes on \emph{Cube 2} with a~$256^{3}$ coefficient matrix; (ii) about 1 hour and 13 minutes with 1 process and about 1 minute on \emph{Cube 3} with a~$512^{3}$ coefficient matrix. For FLOPS, tests outlined 17\% efficiency with \emph{Cube~1}, 42\% with \emph{Cube 2} and 54\% with \emph{Cube 3} when considering 576 processes. This result derives from an improved balancing between computation and communications: Fig.~\ref{fig:granularity} provides the percentage of time spent in MPI communications and computation, while varying the number of processes and the size of the coefficient matrix. Better ratios are achieved with more demanding data sets and a variable number of processes, as reflected in FLOPS efficiency values. 
\subsection{Error analysis\label{sec:ERROR-ANALYSIS}}
To study how the number of iterations and the stopping criteria influence the accuracy and the execution time of iterative solvers, we set the relative error for the stopping condition in Eq. (\ref{eq:relativeError}) equal to~$\epsilon:=10^{-8}, 10^{-12}, 10^{-15}$. For direct solvers based on the SuperLU, we consider the error metrics (\ref{eq:ERROR-METRIC}) with \mbox{$x_{err}\approx 10^{-15}$}. Our tests are performed with 72 parallel processes, i.e., 2 node of Marconi Supercomputer and the \emph{Cube~1} domain.

\begin{figure}[t]
\centering
\includegraphics[height=170pt]{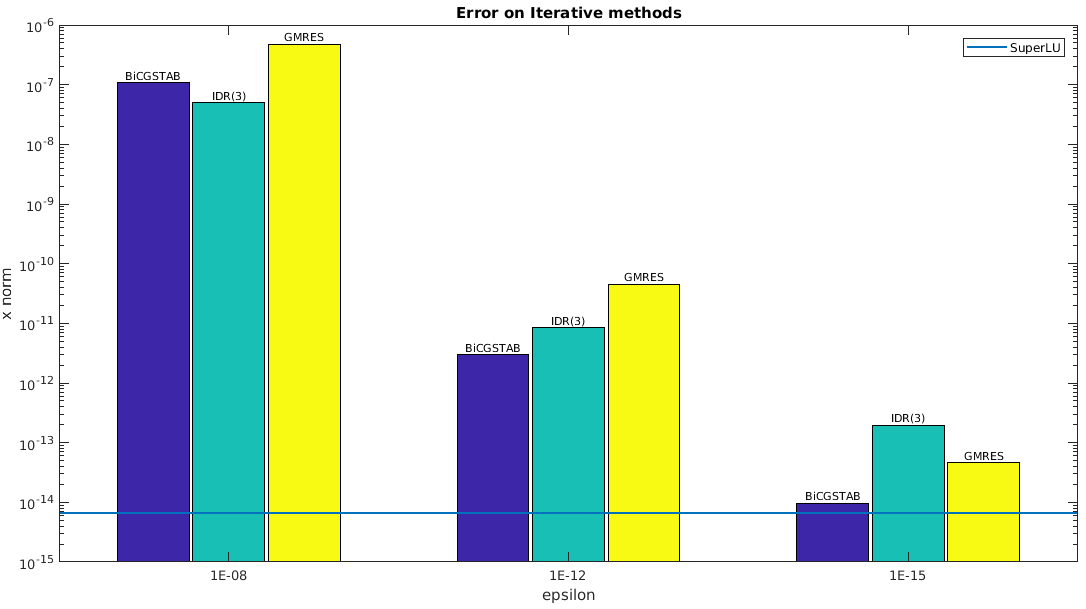}
\caption{Solution accuracy ($y-$axis) of different  iterative  solvers without preconditioning, with respect to different values of~$\epsilon$ ($x$-axis).\label{fig:errorIterativeMethods}}
\end{figure}

Fig.~\ref{fig:errorIterativeMethods} compares the approximation accuracy of different iterative solvers without preconditioning; in particular, we consider BICGSTAB, GMRES, and IDR solvers. This comparison is intended to understand if the selection of a specific iterative algorithm may have an impact on the accuracy. Here, the horizontal line represents the approximation error~$x_{err}$ of SuperLU. Iterative solvers may differ of one order in the approximation error and BICGSTAB generally shows better results, which are further investigated by considering additional preconditioners. 

\begin{figure}[t]
\centering
\includegraphics[height=200pt]{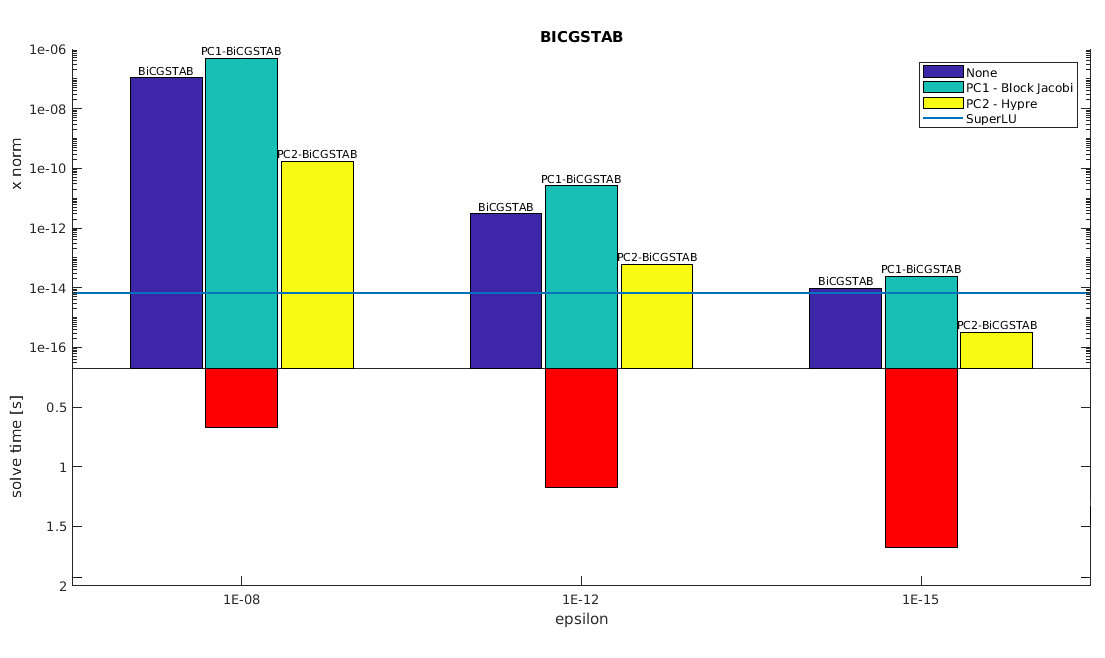}
\caption{($y-$axis) Solution accuracy and timing of the BICGSTAB solver with different preconditioners and values of~$\epsilon$ ($x$-axis).\label{fig:Error}}
\end{figure}

Fig.~\ref{fig:Error} compares the approximation accuracy of BICGSTAB with the Block-Jacobi and Hypre preconditioners and without preconditioners, by varying~$\epsilon$. Hypre needs a larger time (about 25 seconds) to reach a very accurate solution, due to the strong preconditioning phase, while Block-Jacobi has a good approximation accuracy against a negligible execution time (0.01 seconds), which leads to a~$50\%$ reduction of the number of iterations with respect to the solution without preconditioning; for sake of readability in the Figure we depict only the latter. The red bars show the solving time of \emph{Block Jacobi - BICGSTAB} for different~$\epsilon$; varying~$\epsilon$ from~$10^{-8}$ to~$10^{-15}$, the error reaches the order of direct solvers paid with a little impact on execution time, which is still not comparable with respect to direct solvers and provides advantages from memory point of view. We further investigate Block Jacobi - BICGSTAB with respect to the approximation accuracy.
\paragraph*{Further options: polynomial basis}
An important parameter in the modelling of the applicative problem is the degree of the polynomial basis used for the FEM discretisation of the input PDE, which affects the sparsity pattern of the coefficient matrix~$\mathbf{A}$, and thus the solution of the sparse linear system. We analyse the impact on the approximation accuracy when increasing the degree of the polynomial basis from degree 1, i.e. the value used in all tests of this work, to degree 4. Table~\ref{tab:basis} reports the impact for each degree, the coefficient matrix~$\mathbf{A}$  size and its non-zero numbers increase when passing from polynomial degree 1 to degree 4, together with an increase of the number of iterations and execution time ($\epsilon$ fixed to~$10^{-12}$). The approximation accuracy improves significantly from degree 1 to 2 slightly affecting execution time, while approximation accuracy remains almost unchanged with degree 3 and 4 with a notable impact on execution time. 
\begin{table}[t]
\centering
\caption{Block Jacobi - BICGSTAB results with different polynomial degrees, in terms of execution time and accuracy ($\epsilon=10^{-12}$).\label{tab:basis}}
{
\resizebox{\textwidth}{!}{
\begin{tabular}{|c|c|c|c|c|c|c|}
\hline
\textbf{Basis} &\textbf{Nodes} &\textbf{Elem.} &\textbf{Non-zeros} &\textbf{Solver} [$s$] &\textbf{Iter.} &\textbf{Sol. accur.} \\ \hline
P1 &86\,425  &86\,425 &602\,969  &0.11 &358 &~$\num{5.0e-6}$\\ \hline
P2 &86\,425 &344\,697 &3\,894\,474 &0.52 & 689 &~$\num{2.3e-10}$ \\ \hline
P3 &86\,425 &774\,817 &13\,148\,592 &3.32 &938 &~$\num{6.0e-10}$ \\ \hline
P4 &86\,425 &1\,376\,785 &32\,319\,425 &10.21  & 1260 &~$\num{1.7e-9}$ \\ \hline
\end{tabular}}
}
\end{table}
The choice of the polynomial degree influences the approximation accuracy of the solution; it could be selected on the basis of application requirements but it is necessary to take into account the increased computational cost.
\subsection{Handling multiple right-hand side terms\label{sec:MULTIPLE-RHS-TERMS}}
Different physical phenomenon can be described through PDEs discretised by systems with the same coefficient matrix and multiple right-hand side (r.h.s.) terms. For instance, we can consider the approximation of the equilibrium configuration of an isotropic membrane, whose boundary changes its geometry with respect to time. Therefore, we evaluate the execution time of direct and iterative solvers of large sparse linear systems with multiple right-hand side terms. Direct solvers become valuable in this case, since most of the operations (e.g., LU factorisation and column permutation of the coefficient matrix) are performed only once while the solution is computed for each r.h.s. term with a backward and forward substitution in linear time. On the contrary, iterative methods solve a new linear system for every right-hand side term, without preserving any computational phase although the coefficient matrix does not change. During the iterations, only the preconditioner is preserved; according to Sect.~\ref{sec:ITERATIVE-METHODS}, the negligible or predominant computation of the preconditioner depends on the selection of a weak or a strong preconditioner. SuperLU provides the following options %
\begin{equation}\label{eq:MRH-TERMS}
\left\{
\begin{array}{ll}
\mathbf{A}\mathbf{x}_{i}=\mathbf{b}_{i}, \,i=1,\ldots,t &\textrm{(single r.h.s terms)};\\ \\
\mathbf{A}\mathbf{X}=\mathbf{B},\quad
\mathbf{X}:=[\mathbf{x}_{1},\ldots,\mathbf{x}_{t}],\quad
\mathbf{B}:=[\mathbf{b}_{1},\dots,\mathbf{b}_{t}] &\textrm{(block of r.h.s terms)}.
\end{array}
\right.
\end{equation}
The first equation corresponds to the case where r.h.s. terms is passed as one vector per time; the second equation supposes that all the r.h.s. vectors are known a-priori (e.g., they do not depend on the solution computed for the previous r.h.s. terms) and that are passed as a \emph{block} of r.h.s. terms; i.e., a matrix~$\mathbf{B}$. Generally, the matrix~$\mathbf{B}$ is dense and non symmetric; therefore, memory limits may be experienced. 

For iterative solvers, we choose~$\epsilon=10^{-8}, 10^{-12}, 10^{-15}$ on regular grids and we choose only~$\epsilon=10^{-8}$ on irregular grids. Tests have been executed considering 2 nodes of Marconi, i.e. 72 parallel processes, on the data-sets \emph{Cube~1} and \emph{Sphere~1}. As performance metric, we propose execution time varying the number of r.h.s. terms; this means that: (i) the lines with a slower slope correspond to faster solvers, (ii) possible intersections between two lines correspond to a change of behaviour, i.e. the number of r.h.s. terms where one solver becomes more efficient (faster) than the other. 

\begin{figure}[t]
\centering
\includegraphics[width=0.9\textwidth]{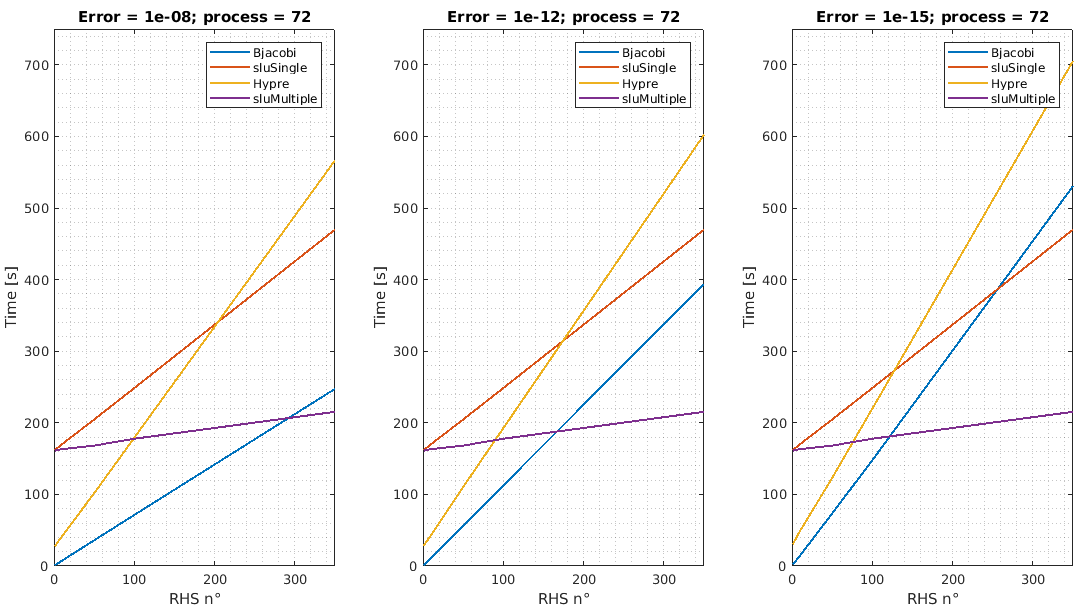}
\caption{Comparison among several approaches in terms of execution time ($y-$axis) with a different number of r.h.s. terms ($x-$axis) on regular grids. \label{fig:multipleRHS}}
\end{figure}

\begin{figure}[t]
\centering
\includegraphics[width=0.90\textwidth]{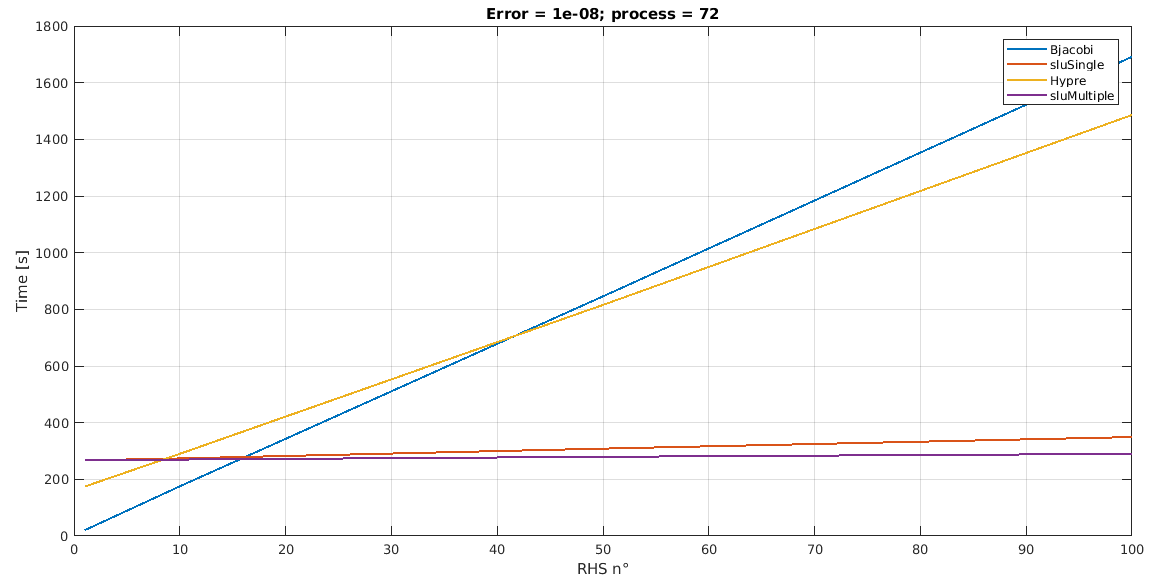}
\caption{Comparison among several approaches in terms of time ($y-$axis) with a different number of r.h.s. terms ($x-$axis) on irregular grids.\label{fig:FFmultipleRHS}}
\end{figure}

We also estimate the number of terms such that direct and iterative solvers require the same computation time. If~$t$ is the number of right-hand side terms, the computation time~$T_{SLU}(t)$ of direct solvers of~$t$ linear systems is equal to the factorisation time~$T_{SLU\_FACT}$ of the coefficient matrix~$\mathbf{A}$ plus the solving time~$T_{SLU\_SOLVE}$ of the two triangular systems for each r.h.s. term; i.e.,
\begin{equation*}
T_{SLU}(t)\approx (T_{SLU\_FACT} + t\cdot T_{SLU\_SOLVE}).
\end{equation*}
Similarly, the computation time~$T_{IS}(t)$ of iterative solvers of~$t$ linear systems is equal to the preconditioning time~$T_{IS\_PC}$ plus the coefficient matrix and the solving time~$T_{IS\_SOLVE}$ for each r.h.s. term; i.e.,
\begin{equation*}
T_{IS}(t)\approx (T_{IS\_PC} + t \cdot T_{IS\_SOLVE}).
\end{equation*}
Then, the number~$t_{0}$ of r.h.s. terms such that the execution time of linear and iterative solvers is comparable is equal to
\begin{equation*}
t_{0}\approx \frac{T_{SLU\_FACT} - T_{IS\_PC}}{T_{IS\_SOLVE} - T_{SLU\_SOLVE}}. 
\end{equation*}
Table~\ref{tab:rsh} reports the value~$t_{0}$ for tests performed on regular/irregular grids and by varying~$\epsilon$. Figs.~\ref{fig:multipleRHS},~\ref{fig:FFmultipleRHS} depict results on regular/irregular grids from direct (two options) and iterative (two preconditioned) solvers. As for regular grids, direct solvers applied to linear systems with multiple r.h.s. terms outperform direct solvers applied to a r.h.s vector per time. BICGSTAB preconditioned with Hypre has the worst performance, since the coefficient matrix is already regular and well distributed. If the exit condition (i.e., the parameter~$\epsilon$) is relaxed, then iterative solvers have a good performance on regular grids and these results are better or comparable with respect to direct solvers; when~$\epsilon$ is reduced, the situation changes. In fact, with~$\epsilon=10^{-8}$, BICGSTAB preconditioned with Block Jacobi (blue line) has a smaller slope than direct method with single r.h.s. terms (red line); direct solvers with a block of r.h.s. terms (violet line) are faster considering more than 290 r.h.s. terms. BICGSTAB preconditioned with Hypre (yellow line) has a larger execution time; indeed, direct solvers obtain more effective results after 210 and 100 r.h.s. terms for the single and block r.h.s. terms, respectively. 
\begin{table}[t]\footnotesize
\centering
\caption{Average number of r.h.s. terms such that iterative and direct solvers have a comparable execution time; Inf and Sup indicate that the method on the horizontal line has always the best/worst performance with respect to vertical intersection.\label{tab:rsh}}
{
\resizebox{\textwidth}{!}{
\begin{tabular}{
p{2.1cm}|
p{0.35cm}p{0.35cm}p{0.35cm}p{0.35cm}|
p{0.35cm}p{0.35cm}p{0.35cm}p{0.35cm}|
p{0.35cm}p{0.35cm}p{0.35cm}p{0.35cm}|
p{0.35cm}p{0.35cm}p{0.35cm}p{0.35cm}|
} 
 & \multicolumn{12}{c|}{Regular grid}                                                                                                                                                                                                                                      & \multicolumn{4}{c|}{Irregular grid}                                                                 \\
 & \multicolumn{4}{c|}{$\epsilon=1e-8$}                                                                              & \multicolumn{4}{c|}{$\epsilon=1e-12$} & \multicolumn{4}{c|}{$\epsilon=1e-15$}                                                         & \multicolumn{4}{c|}{$\epsilon=1e-8$}                                                                             \\ \cline{2-17} 
 & \rotatebox{90}{\textbf{BJacobi}} & \rotatebox{90}{\textbf{Hypre}} & \rotatebox{90}{\textbf{SLU-single}} & \rotatebox{90}{\textbf{SLU-multiple}} & \rotatebox{90}{\textbf{BJacobi}} & \rotatebox{90}{\textbf{Hypre}} & \rotatebox{90}{\textbf{SLU-single}} & \rotatebox{90}{\textbf{SLU-multiple}} & \rotatebox{90}{\textbf{BJacobi}} & \rotatebox{90}{\textbf{Hypre}} & \rotatebox{90}{\textbf{SLU-single}} & \rotatebox{90}{\textbf{SLU-multiple}} & \rotatebox{90}{\textbf{BJacobi}} & \rotatebox{90}{\textbf{Hypre}} & \rotatebox{90}{\textbf{SLU-single}} & \rotatebox{90}{\textbf{SLU-multiple}} \\ \hline 
  \textbf{BJacobi} & - & inf & inf  & 290 & - & inf & 680 & 165 & - & inf & 260 & 120 & - & 41 & 16 & 15 \\  \hline
 \textbf{Hypre} & sup & - & 210 & 100 & sup & - & 175 & 90 & sup & - & 125 & 75 & 41 & - & 9 & 8 \\  \hline
\textbf{SLU-single} & sup  & 210 & - & sup & 680 & 175 & - & sup & 260 & 125 & - & sup & 16 & 9 & - & sup \\ \hline
\textbf{SLU-multiple} & 290  & 100 & inf &- & 165 & 90 & inf & - & 120 & 75 & inf & - & 15 & 8 & inf & -
\\ \hline
\end{tabular}
}
}
\end{table}

Increasing the solution accuracy ($\epsilon=10^{-15}$), iterative solvers become less competitive; direct method with a block of r.h.s. terms provides the best results with respect to the BICGSTAB method preconditioned with Block Jacobi and with Hypre with more than 120 and 75 r.h.s. terms, respectively. BICGSTAB preconditioned with Block Jacobi performs better than direct solvers with single r.h.s. terms up to 260 vectors, while BICGSTAB preconditioned with Hypre is interesting up to 125 r.h.s. terms. Intermediate~$N$ values are obtained for~$\epsilon=10^{-12}$; results are reported in Table~\ref{tab:rsh}. 

Our tests (Fig.~\ref{fig:FFmultipleRHS}) on an irregularly tessellated domain and with~$\epsilon=10^{-8}$ lead to different results. Iterative solvers result less competitive even with low~$\epsilon$ values on irregular grids. BICGSTAB preconditioned with Block Jacobi (blue line) and with Hypre (yellow line) are effective only with a very low number r.h.s. terms; i.e. 16 and 8, respectively. This result is due to their poor scalability figures and to the very slow convergence on irregular grids, as already described in Sect.~\ref{sec:ITERATIVE-METHODS}. If we select an iterative solver, then a strong preconditioner (i.e., Hypre) provides good results in terms of execution time when more than 40 r.h.s. terms are selected. On the basis of these results, we do not further decrease the value of ~$\epsilon$ on irregular domain. 

\section{Discussion} \label{sec:DISCUSSION}
\paragraph*{Direct versus iterative solvers}
SuperLU has good efficiency results on both regular and irregular grids, with an improvement on the latter. On regular grids, SuperLU has a shorter computing timing due to a lower number of non-zero entries. On irregular grids, SuperLU has better scalability results; in fact, a higher non-zero density improves the arithmetic density. The irregular sparsity pattern is managed by a column reordering, which allows a very good distribution of the data. The actual bottleneck of SuperLU is represented by memory requirements that are due to the LU factorisation; therefore, architectural limits can be experienced. This situation may limit the exploitation of SuperLU in terms of the size of the problem/application that can be solved, and in terms of memory capacity of the underlying computational resources. However,  SuperLU is particularly efficient for the solution of large sparse linear systems with multiple r.h.s. terms, and the approximation accuracy is close to the round-off error, when considering well-conditioned coefficient matrices.

As for iterative solvers, better performances are obtained on regular grids, while irregular grids experience an increased number of iterations, due to a more irregular sparsity pattern of the coefficient matrix. We observe that the main operations of the solver do not scale when increasing the number of processes; this aspect has an impact on the solving time and, consequently, on the algorithm efficiency. Beside that, execution time of iterative solvers is considerably lower with respect to direct method ones and memory limits never occurred during our tests. The potential drawback of iterative solvers is a generally lower numerical accuracy, which depends on the number of iterations and the threshold selected for the stopping condition in Eq. (\ref{eq:relativeError}). This situation can be limited by the means of a strong preconditioner and/or a higher polynomial basis degree in the discretisation of Laplace equation, both paid with an increased computational cost and execution time. 
\begin{sidewaystable}[t]\footnotesize
\centering
\caption{Main characteristics of direct and iterative linear solvers.\label{tab:summary}}
\begin{tabular}{@{\extracolsep{4pt}}p{1.3cm}p{3cm}p{3cm}p{3cm}p{1.4cm}p{0.3cm}p{0.3cm}p{0.3cm}p{0.3cm}p{0.3cm}@{}}
\quad & \textbf{Accuracy} & \textbf{Execution time} & \textbf{\makecell[l]{Memory \\requirements} } & \textbf{Data set} & \multicolumn{2}{l}{\textbf{\makecell[l]{Number \\r.h.s. terms}} } & \multicolumn{3}{l}{\textbf{Architecture} }   \\
\cline{6-7} \cline{8-10} 
\multicolumn{1}{l}{\quad} & \multicolumn{1}{l}{\quad} & \multicolumn{1}{l}{\quad} & \multicolumn{1}{l}{\quad} & \multicolumn{1}{l}{\quad} & \rotatebox{45}{One} & \rotatebox{45}{More} & \rotatebox{45}{Embedded} & \rotatebox{45}{Medium}  &  \rotatebox{45}{Large} \\ [2ex]
\hline
\quad & \quad &  \quad & \quad  &  \quad   & \quad & \quad  &  \quad &\quad  & \quad  \\
\textbf{Direct Solver}
&High numerical accuracy \newline 
The approximation accuracy is close to the round-off error, in case of well-conditioned or preconditioned coefficient matrices. &Expensive LU factorisation \newline 
Fast linear solver \newline \newline Execution time is competitive with respect to other direct or iterative solvers in case of multiple r.h.s. terms. \newline \newline 
High scalability on regular and irregular grids &High memory overload for the storage of the~$(\mathbf{L},\mathbf{U})$ factors.
&Regular Grid \newline \newline \newline \newline 
Irregular Grid &\quad & \checkmark \newline \newline \newline \newline \newline \checkmark &  \quad & \checkmark \newline \newline \newline \newline \newline \checkmark & \checkmark \newline \newline \newline \newline \newline \checkmark \\[6ex]
\quad & \quad &   \quad & \quad  &  \quad   & \quad & \quad  &  \quad &\quad  & \quad  \\
\hline
\quad & \quad &   \quad & \quad  &  \quad   & \quad & \quad  &  \quad &\quad  & \quad  \\
\textbf{Iterative Solver} &  Low numerical accuracy but tuneable according to stopping criteria (e.g., number of iterations, numerical accuracy).\newline 
& 
Low numerical accuracy \newline  The method is valuable when solving a single system. \newline \newline
In case of~$t$ r.h.s. terms, the computational cost increases linearly with respect to~$t$. \newline \newline
Higher execution time in case of irregular grids. & Low memory overload, once the coefficient matrix has been stored. Best choice in case of very large sparse linear systems, whose LU factorisation does not fit into the main memory.& 
Regular Grid \newline \newline \newline \newline  Irregular Grid \newline \newline &
\checkmark \newline \newline \newline \newline \newline \checkmark &  \quad & \checkmark \newline \newline \newline \newline \newline \checkmark & \checkmark \newline \newline \newline \newline \newline \checkmark & \checkmark  \\
\quad & \quad &   \quad & \quad  &  \quad   & \quad & \quad  &  \quad &\quad  & \quad  \\
\hline
\end{tabular}
\end{sidewaystable}
\paragraph*{Moving towards multiple r.h.s. terms}
In case of multiple r.h.s. terms, direct solvers generally outperform iterative solvers even when the threshold on accuracy is relatively low. Indeed, iterative solvers are more effective than direct solvers, in case of a low number of r.s.h. terms. If the number of r.s.h. terms grows, direct solvers are more efficient than iterative solvers. In Sect.~\ref{sec:MULTIPLE-RHS-TERMS}, we have identified the average number of terms that make the computation time of direct and iterative solvers comparable; in this case, we have considered each operation of both methods (i.e., factorisation, preconditioning) as a black box, by measuring the execution time without making any assumption on the performed operations. 

We can further analyse this aspect from a computational point of view, i.e., by identifying the average number of r.h.s. terms such that direct solvers have a computational cost comparable with respect to iterative solvers. In case of~$t$ r.h.s. terms and an input coefficient matrix with~$m$ non-zero entries, we have that direct solvers take \mbox{$\mathcal{O}(m)$}-time for the LU factorisation of the input coefficient matrix, \mbox{$\mathcal{O}(n)$}-time for the solution of a linear system with one r.h.s. term through the backward and forward substitution. Indeed, the overall computational cost is \mbox{$\mathcal{C}_{ds}:=\mathcal{O}(m+tn)$}. Iterative solvers solve one linear systems for each r.h.s. term; indeed, the overall computational cost is \mbox{$\mathcal{O}(\sum_{i=1}^{t}k(i)n)$}, where \mbox{$k(i)$} is the number of iterations for the solution of the linear system with the~$i$-th r.h.s. term. Indicating with~$k_{max}$ the maximum number of iterations, the computational cost is \mbox{$\mathcal{C}_{is}:=\mathcal{O}(k_{\max}tn)$}.

The computational cost of direct solvers is lower than iterative solvers if \mbox{$\mathcal{C}_{is}>\mathcal{C}_{ds}$}, i.e., \mbox{$t> t_{0}:=m/(n(k_{\max}-1))$}. Otherwise, direct solvers are computationally more expensive than iterative solvers. If the iterative solver is slowly converging (i.e.,~$k_{\max}$ is large) and the LU factorisation fits into the main memory, then direct solvers are generally more efficient. This remark and the results of our experiments (Sect.~\ref{sec:EXPERIMENT-RESULTS}) confirm the importance of the preconditioning of the coefficient matrix for iterative solvers, which will affect the value~$t_{0}$ of a factor \mbox{$\varphi(n)/(k_{max}-1)$}, where \mbox{$\varphi(n)$} is the cost of the preconditioning step (e.g., \mbox{$\varphi(n)=\mathcal{O}(n)$}, for most of the preconditioners).

\paragraph*{Tools usability}
Another interesting aspect is the usability of the tools in terms of developer' skills and programming effort. SuperLU and PETSc manage the parallelism in a transparent manner with respect to the users, and provide statistics for a high-level profiling of the performances. Since both tools offer the possibility to solve sparse linear systems with few lines of code, their easiness of use is comparable. PETSc represents a more general framework to solve large sparse linear systems and it can be integrated with different external tools, including SuperLU. These considerations imply a larger community of users, a more detailed documentation available online, and support forums. Finally, the performances of SuperLU, used as a stand alone tool or through PETSc interface, are similar. 

\paragraph*{Rule of thumbs}
We summarise these considerations in Table~\ref{tab:summary} (pag.~\pageref{tab:summary}), which presents a general rule of thumbs aimed at supporting application specialists. It provides an empirical base to understand the main peculiarities of direct and iterative solvers, thus driving the effort of application specialists towards the selection of the solver that better fits their requirements. 

Taking into account the application requirements (e.g., number of r.h.s. terms, coefficient matrix characteristics, target approximation accuracy) together with the underlying computational resources (e.g., memory capacity), application specialists have a first insight on the possible options. If the PDEs have to be run on embedded systems (i.e., a resource with limited capabilities), then iterative solvers should be considered as first option. If the input problem requires a high accuracy, then direct solvers should be preferred. If the problem leads to a very large coefficient matrix, then iterative solvers should be considered. In case of multiple r.h.s. terms, previous considerations and in Sect.~\ref{sec:MULTIPLE-RHS-TERMS} guide the user in the selection of the solver with the lowest computation time. 

\section{Conclusions}\label{sec:CONCLUSION}
The solution of PDEs through approximate techniques typically results in the solution of large sparse linear systems. In this context, we selected the Laplace equation on a 2D/3D domain as a standard problem for several physical applications, and we proposed a study about widespread parallel tools for the solution of sparse linear systems with direct and iterative solver. We evaluated their performances by taking into account specific peculiarities of the application domain, e.g. the discretisation of the problem on regular and irregular input domain, the sparsity of the coefficient matrices, etc. We have further compared direct and iterative solvers in terms of memory requirements, numerical accuracy, and their effectiveness when handling of multiple r.h.s. terms. 

The final aim of the work is to support application-domain specialists by providing a first insight on how solvers respond to their applications without the need of performing a large number of experimentation and tests. The choice between direct/iterative solvers has to take into account the knowledge of the application specialists with respect to the underlying problem, and can be based on different elements, such as the properties of the coefficient matrix (e.g., sparsity, symmetry, positive definiteness), single/multiple r.h.s. terms, and matrix size. According to these considerations, we proposed a general rule of thumbs targeted on production-oriented programmers, i.e. application-domain specialists exploiting off-the-shelf tools provided by optimisation-oriented programmers.
Our methodological approach can be applied with different classes of PDEs or with approximation schemes. To the best of our knowledge, a systematic analysis and the definition of a rule of thumbs have not been proposed by scientific community while provides a pragmatic and convenient support in the everyday work.

As future work, we plan to explore the exploitation of GPUs-based tools for the solution of sparse linear systems; we are also interested in improve the understanding of the influence of the discretisation of the input PDE with piecewise linear and polynomial basis functions. In this paper, we mainly relates it with the concept of accuracy, while we would study the wider impact on computational costs and performance. 
\section*{Acknowledgments}
This work has been partially supported by the Programme ``\emph{Advanced Studies in Scientific Computing}'', University La Sapienza - Rome, and the H2020 ERC Advanced Grant CHANGE, Contract 694515. Tests have been supported by CINECA through the ISCRA-C Project HPC-PDE. 


\end{document}